\documentclass[12pt]{article}

\usepackage{amsmath, epsfig, cite}
\usepackage{amssymb}
\usepackage{amsfonts}
\usepackage{latexsym}
\usepackage{amsthm}

\newtheorem{thm}{Theorem}[section]

\numberwithin{equation}{section}

\renewcommand{\thefootnote}

\begin{document}

\begin{center}
{\large\bf On two conjectural series for $\pi$ and their
$q$-analogues
 \footnote{ The work is supported by the National Natural Science Foundation of China (No. 12071103).}}
\end{center}

\renewcommand{\thefootnote}{$\dagger$}

\vskip 2mm \centerline{Chuanan Wei}
\begin{center}
{$^{1}$School of Biomedical Information and Engineering\\ Hainan Medical University, Haikou 571199, China
\\
 Email address: weichuanan78@163.com}
\end{center}

%%date: January 4, 2011
%\vskip 5mm
%\noindent {\it Suggested Running title}: Two Identities of Gould

\vskip 0.7cm \noindent{\bf Abstract.} In terms of the operator
method, we prove two conjectural series for $\pi$ of Sun involving
harmonic numbers of order two. Furthermore, we also give
$q$-analogues of six $\pi$-formulas including the two ones just
mentioned.



\section{Introduction}

For a complex variable $x$, define the well-known Gamma function to
be
\begin{align*}
\Gamma(x)=\int_{0}^{\infty}t^{x-1}e^{-t}dt\quad\text{with}\quad
Re(x)>0.
\end{align*}
Three important properties of it can be stated as follows:
\begin{align*}
\Gamma(x+1)=x\Gamma(x),\quad
\Gamma(x)\Gamma(1-x)=\frac{\pi}{\sin(\pi x)}, \quad\lim_{n\to
\infty}\frac{\Gamma(x+n)}{\Gamma(y+n)}n^{y-x}=1,
\end{align*}
which will be used directly in this paper. For a nonnegative integer
$n$,  define the shifted-factorial as
\begin{align*}
(x)_n=\Gamma(x+n)/\Gamma(x).
\end{align*}
Then we can provide the definition of the hypergeometric series
$$
_{r}F_{s}\left[\begin{array}{c}
a_1,a_2,\ldots,a_{r}\\
b_1,b_2,\ldots,b_{s}
\end{array};\, z
\right] =\sum_{k=0}^{\infty}\frac{(a_1)_k(a_2)_k\cdots(a_{r})_k}
{(b_1)_k(b_2)_k\cdots(b_{s})_k}\frac{z^k}{k!}.
$$

In 1914, Ramanujan \cite{Ramanujan} listed 17 series for $1/\pi$
without proof. Decades later, Borweins \cite{Borwein} proved all of
them firstly. Three of Ramanujan's formulas are expressed as
\begin{align}
&\sum_{k=0}^{\infty}(6k+1)\frac{(\frac{1}{2})_k^3 }{k!^3 4^k}
=\frac{4}{\pi},
 \label{Ramanujan-a}
 \end{align}

 \begin{align}
&\sum_{k=0}^{\infty}(8k+1)\frac{(\frac{1}{2})_k(\frac{1}{4})_k(\frac{3}{4})_k}{k!^3
9^k} =\frac{2\sqrt{3}}{\pi},
 \label{Ramanujan-b}\\
&\quad\:\:\:\sum_{k=0}^{\infty}(42k+5)\frac{(\frac{1}{2})_k^3 }{k!^3
64^k} =\frac{16}{\pi}.
 \label{Ramanujan-c}
\end{align}
There are a lot of different $\pi$-formulas in the literature. Two
of them (cf. \cite[Equation (23)]{Weisstein} and  \cite[P.
221]{Guillera-a}) read
\begin{align}
&\quad\sum_{k=0}^{\infty}\frac{k!}{(2k+1)!!} =\frac{\pi}{2},
 \label{Weisstein}\\
&\sum_{k=0}^{\infty}(3k+2)\frac{(1)_k^3}{(\frac{3}{2})_k^34^k}
=\frac{\pi^2}{4},
 \label{Guillera}
\end{align}
where the double factorial has been defined by
$$\quad (1+2k)!!=\frac{(2k+1)!}{2^kk!}.$$

In 2021, Guo and Lian \cite{Guo-c} conjectured the interesting
double series for $\pi$ related to \eqref{Ramanujan-a}:
\begin{align}
 \sum_{k=0}^{\infty}(6k+1)\frac{(\frac{1}{2})_k^3 }{k!^3
4^k}\sum_{j=1}^{k}\bigg\{\frac{1}{(2j-1)^2}-\frac{1}{16j^2}\bigg\}
=\frac{\pi}{12},
 \label{equation-a}
\end{align}
which has been proved by the auhtor \cite{Wei-a}. Moreover, the
author and Ruan \cite{Wei-b} discovered the following double series
for $\pi$ associated with \eqref{Ramanujan-b}:
\begin{align}
\sum_{k=0}^{\infty}(8k+1)\frac{(\frac{1}{2})_k(\frac{1}{4})_k(\frac{3}{4})_k}{k!^3\,
9^k} \sum_{i=1}^{k}\bigg\{\frac{1}{(2i-1)^2}-\frac{1}{36i^2}\bigg\}
=\frac{\sqrt{3}\,\pi}{54}. \label{equation-b}
\end{align}
 For more
known series on $\pi$, we refer the reader to the papers
\cite{Guillera-b,Schlosser,Chan,Liu,Wang,Zudilin}.

For a complex variable $x$ and two positive integers $\ell,n$,
define the generalized harmonic number of order $\ell$ to be
\[H_{n}^{(\ell)}(x)=\sum_{k=1}^n\frac{1}{(x+k)^{\ell}}.\]
 When $x=0$, it becomes the harmonic number of order $\ell$:
\[H_{n}^{(\ell)}=\sum_{k=1}^n\frac{1}{k^{\ell}}.\]
  Taking $\ell=1$ in $H_{n}^{\langle \ell\rangle}(x)$, we have the generalized harmonic
  number:
\[H_{n}(x)=\sum_{k=1}^n\frac{1}{x+k}.\]  The $x=0$ case of it is the
classical harmonic number:
\[H_{n}=\sum_{k=1}^n\frac{1}{k}.\]
In 2015, Sun \cite{Sun-a} proved a nice series for $\pi^3$
containing harmonic number of order two related to
\eqref{Weisstein}:
\begin{align}
\sum_{k=0}^{\infty}\frac{k!}{(2k+1)!!}H_{k}^{(2)}=\frac{\pi^3}{48}.
 \label{sun-a}
\end{align}
In a recent paper  \cite{Sun-c}, he rewrote \eqref{equation-a} and
\eqref{equation-b} as
\begin{align*}
&\qquad\sum_{k=0}^{\infty}(6k+1)\frac{(\frac{1}{2})_k^3 }{k!^3
4^k}\bigg\{H_{2k}^{(2)}-\frac{5}{16}H_{k}^{(2)}\bigg\}
=\frac{\pi}{12},\\[2mm]
&\sum_{k=0}^{\infty}(8k+1)\frac{(\frac{1}{2})_k(\frac{1}{4})_k(\frac{3}{4})_k}{k!^3\,
9^k} \bigg\{H_{2k}^{(2)}-\frac{5}{18}H_{k}^{(2)}\bigg\}
=\frac{\sqrt{3}\,\pi}{54},
\end{align*}
and proposed the following two conjectures associated with
 \eqref{Guillera} and \eqref{Ramanujan-c} (cf. \cite[Equations (3.67) and
(3.13)]{Sun-c}).

\begin{thm}\label{thm-a}
\begin{align}
\sum_{k=0}^{\infty}(3k+2)\frac{(1)_k^3}{(\frac{3}{2})_k^34^k}\bigg\{H_{2k+1}^{(2)}-\frac{5}{4}H_{k}^{(2)}\bigg\}
=\frac{\pi^4}{48}.
 \label{eq:wei-a}
\end{align}
\end{thm}

\begin{thm}\label{thm-b}
\begin{align}
&\sum_{k=0}^{\infty}(42k+5)\frac{(\frac{1}{2})_k^3 }{k!^3 64^k}
\bigg\{H_{2k}^{(2)}-\frac{25}{92}H_{k}^{(2)}\bigg\}
=\frac{2\pi}{69}.
 \label{eq:wei-b}
\end{align}
\end{thm}

For an integer $n$ and two complex numbers $x$, $q$ with $|q|<1$,
define the $q$-shifted factorial as
\begin{align*}
(x;q)_{\infty}=\prod_{i=0}^{\infty}(1-xq^i),\quad
(x;q)_n=\frac{(x;q)_{\infty}}{(xq^n;q)_{\infty}}.
\end{align*}
For convenience, we sometimes utilize  the compact notation:
\begin{equation*}
(x_1,x_2,\dots,x_r;q)_{m}=(x_1;q)_{m}(x_2;q)_{m}\cdots(x_r;q)_{m},
 \end{equation*}
where $r\in\mathbb{Z}^{+}$ and $m\in\mathbb{Z}^{+}\cup\{0,\infty\}$.
Then following Gasper and Rahman \cite{Gasper}, the basic
hypergeometric series can be defined by
$$
_{r}\phi_{s}\left[\begin{array}{c}
a_1,a_2,\ldots,a_{r}\\
b_1,b_2,\ldots,b_{s}
\end{array};q,\, z
\right] =\sum_{k=0}^{\infty}\frac{(a_1,a_2,\ldots, a_{r};q)_k}
{(q,b_1,b_2,\ldots,b_{s};q)_k}\bigg\{(-1)^kq^{\binom{k}{2}}\bigg\}^{1+s-r}z^k.
$$

Let $[n]=1+q+\cdots+q^{n-1}$ be the $q$-integer. Recently, Guo and
Liu \cite{Guo-d} and Guo and Zudilin \cite{Guo-e} obtained the
following $q$-analogues of \eqref{Ramanujan-a} and
\eqref{Ramanujan-b}:
 \begin{align*}
 &\sum_{k=0}^{\infty}q^{k^2}[6k+1]\frac{(q;q^2)_k^2(q^2;q^4)_k}{(q^4;q^4)_k^3}
=\frac{(1+q)(q^2,q^6;q^4)_{\infty}}{(q^4;q^4)_{\infty}^2},
\\
&\sum_{k=0}^{\infty}q^{2k^2}[8k+1]\frac{(q;q^2)_{k}^2(q;q^2)_{2k}}{(q^2;q^2)_{2k}(q^6;q^6)_{k}^2}
=\frac{(q^3;q^2)_{\infty}(q^3;q^6)_{\infty}}{(q^2;q^2)_{\infty}(q^6;q^6)_{\infty}}.
 \end{align*}
The author \cite{Wei-a} and the author and Ruan \cite{Wei-b} got the
following $q$-analogues of \eqref{equation-a} and
\eqref{equation-b}:
\begin{align*}
&\sum_{k=0}^{\infty}q^{k^2}[6k+1]\frac{(q;q^2)_k^2(q^2;q^4)_k}{(q^4;q^4)_k^3}
\sum_{j=1}^{k}\bigg\{\frac{q^{2j-1}}{[2j-1]^2}-\frac{q^{4j}}{[4j]^2}\bigg\}
\\[3pt]
&\:=\frac{(q^2;q^4)_{\infty}^2(q^5;q^4)_{\infty}}{(q;q^4)_{\infty}(q^4;q^4)_{\infty}^2}
\sum_{i=1}^{\infty}(-1)^{i-1}\frac{q^{2i}}{[2i]^2},
\\[3pt]
&\sum_{k=0}^{\infty}q^{2k^2}[8k+1]\frac{(q;q^2)_{k}^2(q;q^2)_{2k}}{(q^2;q^2)_{2k}(q^6;q^6)_{k}^2}
\sum_{i=1}^{k}\bigg\{\frac{q^{2i-1}}{[2i-1]^2}-\frac{q^{6i}}{[6i]^2}\bigg\}
\notag\\[3pt]
&\:=\frac{(q^3;q^2)_{\infty}(q^3;q^6)_{\infty}}{(q^2;q^2)_{\infty}(q^6;q^6)_{\infty}}
\sum_{j=1}^{\infty}(-1)^{j-1}\frac{q^{3j}}{[3j]^2}.
\end{align*}
More $q$-analogues of
 $\pi$-formulas can be seen in the papers
 \cite{Guo-b,Hou-a,Hou-b,Sun-b}.

 Inspired by the works
 just mentioned, we shall establish $q$-analogues of
\eqref{Weisstein}, \eqref{Guillera}, and \eqref{Ramanujan-c} in the
following theorem.

\begin{thm}\label{thm-c}
\begin{align}
  &\qquad\qquad\qquad\qquad\qquad
  \sum_{k=0}^{\infty}q^{\binom{k+1}{2}}\frac{(q;q)_{k}}{(q^3;q^2)_{k}}
=\frac{(q^2;q^2)_{\infty}^2}{(q,q^3;q^2)_{\infty}},
  \label{eq:wei-d}\\[2mm]
  &\qquad\qquad\qquad
  \sum_{k=0}^{\infty}q^{\binom{k+1}{2}}[3k+2]\frac{(q;q)_{k}^2(q^{2};q^2)_{k}}{(q^3;q^2)_{k}^3}
=\frac{(q^2;q^2)_{\infty}^4}{(q;q^2)_{\infty}^2(q^3;q^2)_{\infty}^2},
  \label{eq:wei-e}\\[2mm]
&\sum_{k=0}^{\infty}q^{6k^2}\frac{(q;q^2)_{k}^3}{(q^2;q^2)_{k}^3}
\frac{(1+q^{1+2k})^3(1-q^{1+6k})-q^{1+6k}(1-q^{3+6k})}{(1-q^3)(1+q)(-q^2;q)_{2k}^3}
=\frac{(q^3,q^5;q^2)_{\infty}}{(q^4;q^2)_{\infty}^2}.
  \label{eq:wei-f}
\end{align}
\end{thm}

Further, we shall furnish $q$-analogues of
\eqref{sun-a}-\eqref{eq:wei-b} in the following three theorems.

\begin{thm}\label{thm-d}
\begin{align}
\sum_{k=0}^{\infty}q^{\binom{k+1}{2}}\frac{(q;q)_{k}}{(q^3;q^2)_{k}}
\sum_{i=1}^{k}\frac{q^{i}}{[i]^2}
=\frac{(q^2;q^2)_{\infty}^2}{(q,q^3;q^2)_{\infty}}
\sum_{j=1}^{\infty}\frac{q^{2j}}{[2j]^2}.
  \label{eq:wei-g}
\end{align}
\end{thm}

\begin{thm}\label{thm-e}
\begin{align}
&\sum_{k=0}^{\infty}q^{\binom{k+1}{2}}[3k+2]\frac{(q;q)_{k}^2(q^{2};q^2)_{k}}{(q^3;q^2)_{k}^3}
\bigg\{\sum_{i=1}^{k}\frac{q^{i}}{[i]^2}-\sum_{i=1}^{k+1}\frac{q^{2i-1}}{[2i-1]^2}\bigg\}
\notag\\[3pt]
&\:=\frac{(q^2;q^2)_{\infty}^4}{(q;q^2)_{\infty}^2(q^3;q^2)_{\infty}^2}
\sum_{j=1}^{\infty}(-1)^{j}\frac{q^{j}}{[j]^2}.
 \label{eq:wei-h}
\end{align}
\end{thm}

\begin{thm}\label{thm-f}
\begin{align}
&\sum_{k=0}^{\infty}q^{6k^2}\frac{(q;q^2)_k^6}{(q^2;q^2)_{2k}^3}\bigg\{\lambda_q(k)\sum_{i=1}^{2k}\frac{q^{2i}}{[2i]^2}
-\mu_q(k)\sum_{i=1}^{k}\frac{q^{2i-1}}{[2i-1]^2}-\nu_q(k)(1-q)q^{1+6k}\bigg\}
\notag\\[2mm]
&\:=\frac{(q,q^3;q^2)_{\infty}}{(q^2;q^2)_{\infty}^2}\bigg\{\sum_{j=1}^{\infty}\frac{q^{2j}}{[2j]^2}-\frac{3(1+q)^3}{64}\sum_{j=1}^{\infty}\frac{q^{2j-1}}{[2j-1]^2}\bigg\},
\label{eq:wei-i}
\end{align}
where
\begin{align*}
&\lambda_q(k)=\frac{1+2q^{1+2k}-q^{1+6k}(2+2q^2+q^{1+2k}+q^{3+2k}-3q^{3+6k})}{(1-q)(1-q^{1+2k})(1+q^{1+2k})^3},
\\[2mm]
&\mu_q(k)=\frac{1+3q^{1+2k}+3q^{2+4k}-2q^{1+6k}+q^{3+6k}-3q^{2+8k}-3q^{3+10k}}{64(1-q)(1+q^{1+2k})^3(137+27q+27q^2+9q^3)^{-1}},
\\[2mm]
&\nu_q(k)=\frac{3(1+q)^3(1+2q^{1+2k}+3q^{2+4k})}{64(1-q^{1+2k})(1+q^{1+2k})^3}-\frac{q^{1+2k}(1+q^{1+2k}+q^{2+4k})^2}{(1-q^{1+2k})(1+q^{1+2k})^5}.
\end{align*}
\end{thm}

For a multivariable function $f(x_1,x_2,\ldots,x_m)$, define the
partial derivative operator $\mathcal{D}_{x_i}$ by
\begin{align*}
&&\mathcal{D}_{x_i}f(x_1,x_2,\ldots,x_m)=\frac{d}{dx_i}f(x_1,x_2,\ldots,x_m)\quad\text{with}\quad
1\leq i\leq m.
 \end{align*}
Then there are the following two relations:
\begin{align*}
&\mathcal{D}_{x}(x+y)_n=(x+y)_nH_n(x+y-1),
\\[3pt]
&\mathcal{D}_{x}(xy;q)_n=-(xy;q)_n\sum_{i=1}^n\frac{yq^{i-1}}{1-xyq^{i-1}}.
\end{align*}

The rest of the paper is arranged as follows. We shall verify
Theorems \ref{thm-a} and \ref{thm-b} via the partial derivative
operator and some summation and transformation formulas for
hypergeometric series in Section 2. Theorems \ref{thm-c}-\ref{thm-f}
will be certified through the partial derivative operator and
several summation and transformation formulas for basic
hypergeometric series in Section 3.

%%%%%%%%%%%%%%%%%%%%%%%%%%%%%%%%%%%%%%%%%%%%%%%%%%%%%%%%%%%%%%%%%%%%%%%%%%%%%%%%%%%%%%%%%%%%%%%%%%%%%%%%%%%%%%%%%%%%%%%%%%%%%%%%%%%%%%%%%%%%%%%%%%%%%%%%%%%%%%%%%%%%%%%%%%%%%%%%%%%%%%%%%%%%%%%%%%%%%%%%%%%%%%%%%%%
\section{Proof of Theorems \ref{thm-a} and \ref{thm-b}}
%%%%%%%%%%%%%%%%%%%%%%%%%%%%%%%%%%%%%%%%%%%%%%%%%%%%%%%%%%%%%%%%%%%%%%%%%%%%%%%%%%%%%%%%%%%%%%%%%%%%%%%%%%%%%%%%%

Above all, we shall prove Theorem \ref{thm-a}.

\begin{proof}[{\bf{Proof of Theorem \ref{thm-a}}}]
 In order to achieve the goal, we need the summation formula for hypergeometric series due to Gosper (1977)(cf. \cite[Equation
 (5.1e)]{Chu}):
\begin{align}
&{_7}F_{6}\left[\begin{array}{c}
a-\frac{1}{2},\frac{2a+2}{3},2b-1,2c-1,2+2a-2b-2c,a+n,-n\\[3pt]
\frac{2a-1}{3},1+a-b,1+a-c,b+c-\frac{1}{2},2a+2n,-2n
\end{array};\, 1 \right]
\notag\\[2mm]
&\quad=\frac{(\frac{1}{2}+a)_n(b)_n(c)_n(a-b-c+\frac{3}{2})_n}
{(\frac{1}{2})_n(1+a-b)_n(1+a-c)_n(b+c-\frac{1}{2})_n}.
\label{equation-c}
\end{align}
 Apply the operator $\mathcal{D}_{b}$ on both sides of the $c=2-b$ case of
\eqref{equation-c} to obtain
\begin{align*}
&\sum_{k=0}^n
\frac{(a-\frac{1}{2})_k(\frac{2a+2}{3})_k(2a-2)_k(2b-1)_k(3-2b)_k(a+n)_k(-n)_k}
{(1)_k(\frac{2a-1}{3})_k(\frac{3}{2})_k(1+a-b)_k(a+b-1)_k(2a+2n)_k(-2n)_k}
\notag\\[2pt]
&\quad\times\Big\{2H_k(2b-2)-2H_k(2-2b)+H_k(a-b)-H_k(a+b-2)\Big\}
\notag\\[2pt]
&=\frac{(a+\frac{1}{2})_n(a-\frac{1}{2})_n(b)_n(2-b)_n}
{(\frac{1}{2})_n(\frac{3}{2})_n(1+a-b)_n(a+b-1)_n}
\notag\\[2pt]
&\quad\times\Big\{H_n(b-1)-H_n(1-b)+H_n(a-b)-H_n(a+b-2)\Big\}.
\end{align*}
Employing the operator $\mathcal{D}_{b}$ on both sides of it, we
have
\begin{align}
&\sum_{k=0}^n
\frac{(a-\frac{1}{2})_k(\frac{2a+2}{3})_k(2a-2)_k(2b-1)_k(3-2b)_k(a+n)_k(-n)_k}
{(1)_k(\frac{2a-1}{3})_k(\frac{3}{2})_k(1+a-b)_k(a+b-1)_k(2a+2n)_k(-2n)_k}
\notag\\[2pt]
&\quad\times\Big\{\Big[2H_k(2b-2)-2H_k(2-2b)+H_k(a-b)-H_k(a+b-2)\Big]^2
\notag\\[2pt]
&\qquad-\Big[4H_k^{(2)}(2b-2)+4H_k^{(2)}(2-2b)-H_k^{(2)}(a-b)-H_k^{(2)}(a+b-2)\Big]\Big\}
\notag\\[2pt]
&=\frac{(a+\frac{1}{2})_n(a-\frac{1}{2})_n(b)_n(2-b)_n}
{(\frac{1}{2})_n(\frac{3}{2})_n(1+a-b)_n(a+b-1)_n}
\notag\\[2pt]
&\quad\times\Big\{\Big[H_n(b-1)-H_n(1-b)+H_n(a-b)-H_n(a+b-2)\Big]^2
\notag\\[2pt]
&\qquad-\Big[H_n^{(2)}(b-1)+H_n^{(2)}(1-b)-H_n^{(2)}(a-b)-H_n^{(2)}(a+b-2)\Big]\Big\}.
\label{eq:wei-ee}
\end{align}
 The $(a,b)=(\frac{3}{2},1)$ case of \eqref{eq:wei-ee} engenders
\begin{align}
&\sum_{k=0}^{n}\frac{(\frac{5}{3})_k(1)_k^3(\frac{3}{2}+n)_k(-n)_k}{(\frac{2}{3})_k(\frac{3}{2})_k^3(3+2n)_k(-2n)_k}
\bigg\{4H_{2k+1}^{(2)}-5H_{k}^{(2)}-4\bigg\}
\notag\\[2pt]
&\:=\frac{\Gamma(2+n)\Gamma(1+n)^3}
{\Gamma(\frac{1}{2}+n)\Gamma(\frac{3}{2}+n)^3}
\frac{\Gamma(\frac{1}{2})\Gamma(\frac{3}{2})^3}{\Gamma(2)\Gamma(1)^3}\bigg\{4H_{2n+1}^{(2)}-2H_{n}^{(2)}-4\bigg\}.
\label{eq:wei-ff}
\end{align}
Since that the $(a,b,c)=(\frac{3}{2},1,1)$ case of
\eqref{equation-c} reads
\begin{align}
\sum_{k=0}^{n}\frac{(\frac{5}{3})_k(1)_k^3(\frac{3}{2}+n)_k(-n)_k}{(\frac{2}{3})_k(\frac{3}{2})_k^3(3+2n)_k(-2n)_k}
=\frac{\Gamma(2+n)\Gamma(1+n)^3}
{\Gamma(\frac{1}{2}+n)\Gamma(\frac{3}{2}+n)^3}
\frac{\Gamma(\frac{1}{2})\Gamma(\frac{3}{2})^3}{\Gamma(2)\Gamma(1)^3},
\label{eq:wei-gg}
\end{align}
the linear combination of \eqref{eq:wei-ff} and \eqref{eq:wei-gg}
gives
\begin{align*}
&\sum_{k=0}^{n}\frac{(\frac{5}{3})_k(1)_k^3(\frac{3}{2}+n)_k(-n)_k}{(\frac{2}{3})_k(\frac{3}{2})_k^3(3+2n)_k(-2n)_k}
\bigg\{4H_{2k+1}^{(2)}-5H_{k}^{(2)}\bigg\}
\notag\\[2pt]
&\:=\frac{\Gamma(2+n)\Gamma(1+n)^3}
{\Gamma(\frac{1}{2}+n)\Gamma(\frac{3}{2}+n)^3}
\frac{\Gamma(\frac{1}{2})\Gamma(\frac{3}{2})^3}{\Gamma(2)\Gamma(1)^3}\bigg\{4H_{2n+1}^{(2)}-2H_{n}^{(2)}\bigg\}.
\end{align*}
Letting $(x,n)\to(\frac{1}{2},\infty)$ and making use of Euer's
formula:
\begin{align}\label{Euler}
\sum_{j=1}^{\infty}\frac{1}{j^2}=\frac{\pi^2}{6},
\end{align}
we catch hold of \eqref{eq:wei-a}.
\end{proof}

Subsequently, we shall display the proof of Theorem \ref{thm-b}.

\begin{proof}[{\bf{Proof of Theorem \ref{thm-b}}}]
Recall a  transformation formula for hypergeometric series (cf.
\cite[Theorem 31]{Chu-b}):
\begin{align*}
&\sum_{k=0}^{\infty}(-1)^k\frac{(b)_k(c)_k(d)_k(e)_k(1+a-b-c)_k(1+a-b-d)_{k}(1+a-b-e)_{k}}{(1+a-b)_{2k}(1+a-c)_{2k}(1+a-d)_{2k}(1+a-e)_{2k}}
\\[2mm]
&\quad\times\frac{(1+a-c-d)_k(1+a-c-e)_{k}(1+a-d-e)_{k}}{(1+2a-b-c-d-e)_{2k}}\sigma_k(a,b,c,d,e)
\\[2mm]
&=\sum_{k=0}^{\infty}(a+2k)\frac{(b)_k(c)_k(d)_k(e)_k}{(1+a-b)_{k}(1+a-c)_{k}(1+a-d)_{k}(1+a-e)_{k}},
\end{align*}
where
\begin{align*}
&\sigma_k(a,b,c,d,e)\\[2mm]
&\quad=\frac{(1+2a-b-c-d+3k)(a-e+2k)}{(1+2a-b-c-d-e+2k)}+\frac{(e+k)(1+a-b-c+k)}{(1+a-b+2k)(1+a-d+2k)}
\\[2mm]
&\qquad\times\frac{(1+a-b-d+k)(1+a-c-d+k)(2+2a-b-d-e+3k)}{(1+2a-b-c-d-e+2k)(2+2a-b-c-d-e+2k)}
\\[2mm]
&\quad+\frac{(c+k)(e+k)(1+a-b-c+k)(1+a-b-d+k)}{(1+a-b+2k)(1+a-c+2k)(1+a-d+2k)(1+a-e+2k)}
\\[2mm]
&\qquad\times\frac{(1+a-b-e+k)(1+a-c-d+k)(1+a-d-e+k)}{(1+2a-b-c-d-e+2k)(2+2a-b-c-d-e+2k)}.
\end{align*}
Choosing $(a,b,c,d,e)=(\frac{1}{2},\frac{1}{2},x,1-x,-n)$ in the
last equation and calculating the series on the right-hand side by
 Dougall's $_5F_4$ summation formula (cf. \cite[P.
71]{Andrews}):
\begin{align*}
{_{5}F_{4}}\left[\begin{array}{cccccccc}
  a,1+\frac{a}{2},b,c,-n\\
  \frac{a}{2},1+a-b,1+a-c,1+a+n
\end{array};1
\right]=\frac{(1+a)_n(1+a-b-c)_n}{(1+a-b)_n(1+a-c)_n},
\end{align*}
we arrive at
\begin{align}
&\sum_{k=0}^{n}(-1)^k\frac{(\frac{1}{2})_k^2(1+n)_k(-n)_k}{(1)_{2k}(\frac{1}{2}+n)_{2k}(\frac{3}{2}+n)_{2k}}
\frac{(x)_k^2(1-x)_k^2(\frac{1}{2}+x+n)_k(\frac{3}{2}-x+n)_k}{(\frac{1}{2}+x)_{2k}(\frac{3}{2}-x)_{2k}}
\notag\\[2mm]
&\quad\times\Omega_k(x;n)=\frac{(\frac{1}{2})_{n}(\frac{3}{2})_{n}}{(\frac{1}{2}+x)_{n}(\frac{3}{2}-x)_{n}},
\label{eq:wei-aa}
\end{align}
where
\begin{align*}
\Omega_k(x;n)
&=(1+6k)+\frac{4(x+k)(1-x+k)(3+2x+2n+6k)(k-n)}{(1+2x+4k)(1+2n+4k)(3+2n+4k)}
\\[2mm]
&\quad+\frac{16(x+k)^2(1-x+k)(1+2x+2n+2k)(1+n+k)(k-n)}{(1+2x+4k)(3-2x+4k)(1+2n+4k)(3+2n+4k)^2}.
\end{align*}
Apply the operator $\mathcal{D}_{x}$ on both sides of
\eqref{eq:wei-aa} to get
\begin{align}
&\sum_{k=0}^{n}(-1)^k\frac{(\frac{1}{2})_k^2(1+n)_k(-n)_k}{(1)_{2k}(\frac{1}{2}+n)_{2k}(\frac{3}{2}+n)_{2k}}
\frac{(x)_k^2(1-x)_k^2(\frac{1}{2}+x+n)_k(\frac{3}{2}-x+n)_k}{(\frac{1}{2}+x)_{2k}(\frac{3}{2}-x)_{2k}}
\notag\\[2mm]
&\quad\times\Big\{2H_k(x-1)-2H_k(-x)-H_{2k}(x-\tfrac{1}{2})+H_{2k}(\tfrac{1}{2}-x)+H_{k}(x-\tfrac{1}{2}+n)
\notag\\[2mm]
&\qquad-H_{k}(\tfrac{1}{2}-x+n)\Big\}\Omega_k(x;n)
\notag\\[2mm]
&\:+\sum_{k=0}^{n}(-1)^k\frac{(\frac{1}{2})_k^2(1+n)_k(-n)_k}{(1)_{2k}(\frac{1}{2}+n)_{2k}(\frac{3}{2}+n)_{2k}}
\frac{(x)_k^2(1-x)_k^2(\frac{1}{2}+x+n)_k(\frac{3}{2}-x+n)_k}{(\frac{1}{2}+x)_{2k}(\frac{3}{2}-x)_{2k}}
\notag\\[2mm]
&\qquad\times\mathcal{D}_{x}\Omega_k(x;n)
=\frac{(\frac{1}{2})_{n}(\frac{3}{2})_{n}}{(\frac{1}{2}+x)_{n}(\frac{3}{2}-x)_{n}}
\Big\{H_{n}(\tfrac{1}{2}-x)-H_{n}(x-\tfrac{1}{2})\Big\}.
\label{eq:wei-bb}
\end{align}
Dividing both sides by $1-2x$ and utilizing the relation
 \begin{align}
 \frac{1}{v-u-2x}\Big\{H_m(x+u)-H_m(v-x)\Big\}=\sum_{i=1}^m\frac{1}{(x+u+i)(v-x+i)},
 \label{relation-a}
\end{align}
Equation \eqref{eq:wei-bb} can be manipulated as

\begin{align*}
&\sum_{k=0}^{n}(-1)^k\frac{(\frac{1}{2})_k^2(1+n)_k(-n)_k}{(1)_{2k}(\frac{1}{2}+n)_{2k}(\frac{3}{2}+n)_{2k}}
\frac{(x)_k^2(1-x)_k^2(\frac{1}{2}+x+n)_k(\frac{3}{2}-x+n)_k}{(\frac{1}{2}+x)_{2k}(\frac{3}{2}-x)_{2k}}
\notag\\[2mm]
&\quad\times\bigg\{2\sum_{i=1}^k\frac{1}{(x-1+i)(-x+i)}-\sum_{i=1}^{2k}\frac{1}{(x-\frac{1}{2}+i)(\frac{1}{2}-x+i)}
\notag\\[2mm]
&\qquad+\sum_{i=1}^k\frac{1}{(x-\frac{1}{2}+n+i)(\frac{1}{2}-x+n+i)}\bigg\}\Omega_k(x;n)
\notag\\[2mm]
&\:+\sum_{k=0}^{n}(-1)^k\frac{(\frac{1}{2})_k^2(1+n)_k(-n)_k}{(1)_{2k}(\frac{1}{2}+n)_{2k}(\frac{3}{2}+n)_{2k}}
\frac{(x)_k^2(1-x)_k^2(\frac{1}{2}+x+n)_k(\frac{3}{2}-x+n)_k}{(\frac{1}{2}+x)_{2k}(\frac{3}{2}-x)_{2k}}
\notag\\[2mm]
&\qquad\times\frac{\mathcal{D}_{x}\Omega_k(x;n)}{1-2x}
=-\frac{(\frac{1}{2})_{n}(\frac{3}{2})_{n}}{(\frac{1}{2}+x)_{n}(\frac{3}{2}-x)_{n}}
\sum_{j=1}^{n}\frac{1}{(x-\frac{1}{2}+j)(\frac{1}{2}-x+j)}.
\end{align*}
Letting $(x,n)\to(\frac{1}{2},\infty)$ and drawing upon Euer's
formula \eqref{Euler}, there is
\begin{align}
&\sum_{k=0}^{\infty}\frac{(\frac{1}{2})_k^3 }{k!^3 64^k}
\bigg\{(42k+5)\Big[2H_{k}^{(2)}-7H_{2k}^{(2)}\Big]+\frac{9}{1+2k}\bigg\}
=\frac{8\pi}{3}.
 \label{eq:wei-cc}
\end{align}

Recollect a summation formula for  hypergeometric series (cf.
\cite[Corollary 2.33]{Chu-a}):
\begin{align}
&\sum_{k=0}^{\infty}\frac{(x)_k^3(1-x)_k^3}{(1)_{k}^3(\frac{3}{2})_{k}^3}
\frac{k(1+3k)(3+9k+7k^2)+x(1-x)(1+6k+6k^2+x-x^2)}{64^k}
\notag\\[2mm]
&\:=\frac{\sin(\pi x)}{\pi}, \label{series-a}
\end{align}
where we have replaced $\sin(\pi x)/\pi x$ by $\sin(\pi x)/\pi $ for
correction. When $0<x<1$, it is obvious that the series on the
left-hand side  is uniformly convergent. Employing the operator
$\mathcal{D}_{x}$ on both sides of of \eqref{series-a} and taking
advantage of \eqref{relation-a}, there holds
\begin{align*}
&3\sum_{k=0}^{\infty}\frac{(x)_k^3(1-x)_k^3}{(1)_{k}^3(\frac{3}{2})_{k}^3}\sum_{i=1}^k\frac{1}{(x-1+i)(-x+i)}
\\[2mm]
&\quad\times\frac{k(1+3k)(3+9k+7k^2)+x(1-x)(1+6k+6k^2+x-x^2)}{64^k}
\\[2mm]
&+\sum_{k=0}^{\infty}\frac{(x)_k^3(1-x)_k^3}{(1)_{k}^3(\frac{3}{2})_{k}^3}\frac{1+6k+6k^2+2x(1-x)}{64^k}
=\frac{\cos(\pi x)}{1-2x}.
\end{align*}
The $x\to\frac{1}{2}$ case of the upper identity provides
\begin{align}
\sum_{k=0}^{\infty}\frac{(\frac{1}{2})_k^3 }{k!^3 64^k}
\bigg\{(42k+5)\Big[4H_{2k}^{(2)}-H_{k}^{(2)}\Big]+\frac{8}{1+2k}\bigg\}
=\frac{8\pi}{3}.
 \label{eq:wei-dd}
\end{align}
Hence we deduce \eqref{eq:wei-b} from the linear combination of
\eqref{eq:wei-cc} and \eqref{eq:wei-dd}.
\end{proof}

\section{Proof of Theorems \ref{thm-c}-\ref{thm-f}}
%%%%%%%%%%%%%%%%%%%%%%%%%%%%%%%%%%%%%%%%%%%%%%%%%%%%%%%%%%%%%%%%%%%%%%%%%%%%%%%%%%%%%%%%%%%%%%%%%%%%%%%%%%%%%%%%%

Firstly, we are ready to prove Theorem \ref{thm-c}.

\begin{proof}[{\bf{Proof of Theorem \ref{thm-c}}}]

For achieving the purpose, we require two identities for basic
hypergeometric series (cf. \cite[Equation (5.1d)]{Chu-q} and
\cite[Theorem 17]{Chu-2012}):
\begin{align}
&\sum_{k=0}^{n}\frac{1-aq^{3k-1}}{1-aq^{-1}}
\frac{(q^{-2n},aq^{2n},a/q;q^2)_k}{(aq^{2}/b,aq^{2}/c,bc/q;q^2)_k}
\frac{(b/q,c/q,aq^{2}/bc;q)_k}{(q,aq^{2n},q^{-2n};q)_k} q^{k}
\notag\\[3pt]
  &\:=\frac{(aq,b,c,aq^{3}/bc;q^2)_n}{(q,aq^{2}/b,aq^{2}/c,bc/q;q^2)_n},
\label{basic}
\end{align}
\begin{align}
&\sum_{k=0}^{\infty}\frac{(b,c,d,e,aq/bc,aq/bd,aq/be,aq/cd,aq/ce,aq/de;q)_{k}}{(aq/b,aq/c,aq/d,aq/e,a^2q/bcde;q)_{2k}}
\notag\\[2pt]
 &\quad\times\:
 q^{\frac{k(5k+1)}{2}}\bigg(\frac{-a^3}{bcde}\bigg)^kA_k(a,b,c,d,e;q)
 \notag\\[2pt]
&\:=\sum_{k=0}^{\infty}\frac{1-aq^{2k}}{1-a}\frac{(b,c,d,e;q)_k}{(aq/b,aq/c,aq/d,aq/e;q)_k}\bigg(\frac{a^2q}{bcde}\bigg)^k,
\label{basic-a}
\end{align}
where
\begin{align*}
&A_k(a,b,c,d,e;q)=\frac{(1-q^{2k}a/e)(1-q^{1+3k}a^2/bcd)}{(1-a)(1-q^{1+2k}a^2/bcde)}+\frac{q^{2k}a}{e}\frac{(1-q^ke)}
{(1-a)}\\[2pt]
&\quad\times\frac{(1-q^{1+k}a/bc)(1-q^{1+k}a/bd)(1-q^{1+k}a/cd)(1-q^{2+3k}a^2/bde)}
{(1-q^{1+2k}a/b)(1-q^{1+2k}a/d)(1-q^{1+2k}a^2/bcde)(1-q^{2+2k}a^2/bcde)}\\[2pt]
&+\frac{q^{1+4k}a^2}{ce}\frac{(1-q^kc)(1-q^ke)(1-q^{1+k}a/bc)(1-q^{1+k}a/bd)(1-q^{1+k}a/be)}{(1-a)(1-q^{1+2k}a/b)(1-q^{1+2k}a/c)(1-q^{1+2k}a/d)(1-q^{1+2k}a/e)}\\
&\quad\times\frac{(1-q^{1+k}a/cd)(1-q^{1+k}a/de)}{(1-q^{1+2k}a^2/bcde)(1-q^{2+2k}a^2/bcde)}.
\end{align*}

Notice that the $(a,b,c)=(0,q^2,q^2)$ case of \eqref{basic} is
\begin{align*}
\sum_{k=0}^{n}\frac{(q;q)_k}{(q^3;q^2)_k}
\frac{(q^{-2n};q^2)_k}{(q^{-2n};q)_k} q^{k}
=\frac{(q^2;q^2)_n^2}{(q,q^3;q^2)_n}.
\end{align*}
Letting $n\to \infty$ in the above identity, we obtain
\eqref{eq:wei-d}.

The $(a,b,c)=(q^3,q^2,q^2)$ case of \eqref{basic} reads
\begin{align*}
\sum_{k=0}^{n}\frac{1-q^{3k+2}}{1-q^{2}}
\frac{(q;q)_k^2(q^2;q^2)_k}{(q^3;q^2)_k^3}
\frac{(q^{2n+3},q^{-2n};q^2)_k}{(q^{2n+3},q^{-2n};q)_k} q^{k}
=\frac{(q^2;q^2)_n^3(q^4;q^2)_n}{(q;q^2)_n(q^3;q^2)_n^3}.
\end{align*}
Letting $n\to \infty$ in the upper identity, we get
\eqref{eq:wei-e}.

Performing the replacements $(a, b, c, d, e)\to(x, a, b,  xq/c,
xq/d)$ in \eqref{basic-a} and then letting $x\to 0$, we find
\begin{align}
&\sum_{k=0}^{\infty}\frac{(a,b,c/a,c/b,d/a,d/b;q)_{k}}{(c,d;q)_{1+2k}(cd/abq;q)_{2+2k}}
 q^{6\binom{k}{2}}\bigg(\frac{c^2d^2}{ab}\bigg)^kB_k(a,b,c,d;q)
 \notag\\[2pt]
&\:=\sum_{k=0}^{\infty}\frac{(a,b;q)_k}{(c,d;q)_k}\bigg(\frac{cd}{abq}\bigg)^k,
\label{basic-c}
\end{align}
where
\begin{align*}
B_k(a,b,c,d;q)&=(1-q^{2k}c)(1-q^{2k}d)(1-q^{2k}cd/ab)(1-q^{3k-1}cd/b)
\\[2pt]
&\quad-\frac{q^{3k-1}cd}{a}(1-q^{k}a)(1-q^{k}c/b)(1-q^{k}d/b)(1-q^{3k}cd/a).
\end{align*}
When $d=q$, the series on the right-hand side of \eqref{basic-c} can
be evaluated by the $q$-Gauss summation formula (cf. \cite[Appendix
II. 8]{Gasper}):
\begin{align*}
{_2}\phi_{1}\left[\begin{array}{c}
a,b\\[3pt]
c
\end{array};\, q, \frac{c}{ab} \right]
=\frac{(c/a,c/b;q)_{\infty}}{(c,c/ab;q)_{\infty}}.
\end{align*}
So we have
\begin{align}
&\sum_{k=0}^{\infty}\frac{(a,b,q/a,q/b,c/a,c/b;q)_{k}}{(q;q)_{1+2k}(cq,cq^2/ab;q)_{2k}}
 q^{3k^2-k}\bigg(\frac{c^2}{ab}\bigg)^kC_k(a,b,c;q)
 \notag\\[2pt]
&\:=\frac{(c/a,c/b;q)_{\infty}}{(cq,cq^2/ab;q)_{\infty}},
\label{basic-d}
\end{align}
where
\begin{align*}
C_k(a,b,c;q)&=(1-q^{2k}c)(1-q^{1+2k})(1-q^{1+2k}c/ab)(1-q^{3k}c/b)
\\[2pt]
&\quad-\frac{q^{3k}c}{a}(1-q^{k}a)(1-q^{k}c/b)(1-q^{1+k}/b)(1-q^{1+3k}c/a).
\end{align*}
Letting $(a, b, c, q)\to(q, q, q^2, q^2)$ in \eqref{basic-d}, we
prove \eqref{eq:wei-f}.

\end{proof}

Secondly, we start to prove Theorem \ref{thm-d}.

\begin{proof}[{\bf{Proof of Theorem \ref{thm-d}}}]
Apply the operator $\mathcal{D}_{b}$ on both sides of the $c\to
q^{4}/b$ case of \eqref{basic} to discover
\begin{align*}
&\sum_{k=0}^{n}\frac{1-aq^{3k-1}}{1-aq^{-1}}
\frac{(q^{-2n},aq^{2n},a/q;q^2)_k}{(q^{3},aq^{2}/b,ab/q^{2};q^2)_k}
\frac{(a/q^2,b/q,q^{3}/b;q)_k}{(q,aq^{2n},q^{-2n};q)_k}
q^{k}D_k(a,b)
\notag\\[3pt]
  &\:=\frac{(aq,a/q,b,q^{4}/b;q^2)_n}{(q,q^{3},aq^{2}/b,ab/q^{2};q^2)_n}E_n(a,b),
\end{align*}
where
\begin{align*}
&D_k(a,b)=\sum_{i=1}^k\frac{q^{i-2}}{1-bq^{i-2}}-\sum_{i=1}^k\frac{q^{i+2}/b^2}{1-q^{i+2}/b}
+\sum_{i=1}^k\frac{aq^{2i}/b^2}{1-aq^{2i}/b}-\sum_{i=1}^k\frac{aq^{2i-4}}{1-abq^{2i-4}},
\\[2pt]
&E_n(a,b)=\sum_{j=1}^n\frac{q^{2j-2}}{1-bq^{2j-2}}-\sum_{j=1}^n\frac{q^{2j+2}/b^2}{1-q^{2j+2}/b}
+\sum_{j=1}^n\frac{aq^{2j}/b^2}{1-aq^{2j}/b}-\sum_{j=1}^n\frac{aq^{2j-4}}{1-abq^{2j-4}}.
\end{align*}
Employing the operator $\mathcal{D}_{b}$ on both sides of the last
equation, it is easy to show that
\begin{align}
&\sum_{k=0}^{n}\frac{1-aq^{3k-1}}{1-aq^{-1}}
\frac{(q^{-2n},aq^{2n},a/q;q^2)_k}{(q^{3},aq^{2}/b,ab/q^{2};q^2)_k}
\frac{(a/q^2,b/q,q^{3}/b;q)_k}{(q,aq^{2n},q^{-2n};q)_k} q^{k}
\bigg\{D_k(a,b)^2-F_k(a,b)\bigg\}
\notag\\[3pt]
&\quad=\frac{(aq,a/q,b,q^{4}/b;q^2)_n}{(q,q^{3},aq^{2}/b,ab/q^{2};q^2)_n}\bigg\{E_n(a,b)^2-G_n(a,b)\bigg\},
\label{eq:wei-eee}
\end{align}
where
\begin{align*}
&F_k(a,b)=\sum_{i=1}^k\frac{q^{2i-4}}{(1-bq^{i-2})^2}-\sum_{i=1}^k\frac{(q^{i+2}/b-2)q^{i+2}/b^3}{(1-q^{i+2}/b)^2}
\\[2pt]
&\qquad\qquad\:+\sum_{i=1}^k\frac{(aq^{2i}/b-2)aq^{2i}/b^3}{(1-aq^{2i}/b)^2}-\sum_{i=1}^k\frac{a^2q^{4i-8}}{(1-abq^{2i-4})^2},
\\[2pt]
&G_n(a,b)=\sum_{j=1}^n\frac{q^{4j-4}}{(1-bq^{2j-2})^2}-\sum_{j=1}^n\frac{(q^{2j+2}/b-2)q^{2j+2}/b^3}{(1-q^{2j+2}/b)^2}
\\[2pt]
&\qquad\qquad\:+\sum_{j=1}^n\frac{(aq^{2j}/b-2)aq^{2j}/b^3}{(1-aq^{2j}/b)^2}-\sum_{j=1}^n\frac{a^2q^{4j-8}}{(1-abq^{2j-4})^2}.
\end{align*}
The $(a,b)=(0,q^2)$ case of \eqref{eq:wei-eee} produces
\begin{align*}
\sum_{k=0}^{n}q^{k}\frac{(q;q)_{k}(q^{-2n};q^2)_{k}}{(q^3;q^2)_{k}(q^{-2n};q)_{k}}
\sum_{i=1}^{k}\frac{q^{i}}{[i]^2}
=\frac{(q^2;q^2)_{n}^2}{(q,q^3;q^2)_{n}}
\sum_{j=1}^{n}\frac{q^{2j}}{[2j]^2}.
\end{align*}
Letting $n\to \infty$ in this identity, we catch hold of
\eqref{eq:wei-g}.
\end{proof}

Thirdly, we shall prove Theorem \ref{thm-e}.

\begin{proof}[{\bf{Proof of Theorem \ref{thm-e}}}]
The $(a, b)=(q^3,q^2)$ case of \eqref{eq:wei-eee} provides
\begin{align}
&\sum_{k=0}^{n}q^{k}\frac{1-q^{3k+2}}{1-q^2}\frac{(q;q)_{k}^2(q^{2};q^2)_{k}}{(q^3;q^2)_{k}^3}
\frac{(q^{2n+3},q^{-2n};q^2)_{k}}{(q^{2n+3},q^{-2n};q)_{k}}\bigg\{\sum_{i=1}^{k}\frac{q^{i}}{[i]^2}-\sum_{i=2}^{k+1}\frac{q^{2i-1}}{[2i-1]^2}\bigg\}
\notag\\[3pt]
&\:=\frac{(q^2;q^2)_{n}^3(q^4;q^2)_{n}}{(q;q^2)_{n}(q^3;q^2)_{n}^3}
\bigg\{\sum_{i=1}^{n}\frac{q^{2j}}{[2j]^2}-\sum_{j=2}^{n+1}\frac{q^{2j-1}}{[2j-1]^2}\bigg\}.
 \label{eq:wei-fff}
\end{align}
The $(a,b,c)=(q^3,q^2,q^2)$ case of \eqref{basic} can be expressed
as
\begin{align}
\sum_{k=0}^{n}q^{k}\frac{1-q^{3k+2}}{1-q^2}\frac{(q;q)_{k}^2(q^{2};q^2)_{k}}{(q^3;q^2)_{k}^3}
\frac{(q^{2n+3},q^{-2n};q^2)_{k}}{(q^{2n+3},q^{-2n};q)_{k}}
=\frac{(q^2;q^2)_{n}^3(q^4;q^2)_{n}}{(q;q^2)_{n}(q^3;q^2)_{n}^3}.
 \label{eq:wei-ggg}
\end{align}
According to the linear combination of \eqref{eq:wei-fff} and
\eqref{eq:wei-ggg}, we have
\begin{align*}
&\sum_{k=0}^{n}q^{k}\frac{1-q^{3k+2}}{1-q^2}\frac{(q;q)_{k}^2(q^{2};q^2)_{k}}{(q^3;q^2)_{k}^3}
\frac{(q^{2n+3},q^{-2n};q^2)_{k}}{(q^{2n+3},q^{-2n};q)_{k}}\bigg\{\sum_{i=1}^{k}\frac{q^{i}}{[i]^2}-\sum_{i=1}^{k+1}\frac{q^{2i-1}}{[2i-1]^2}\bigg\}
\notag\\[3pt]
&\:=\frac{(q^2;q^2)_{n}^3(q^4;q^2)_{n}}{(q;q^2)_{n}(q^3;q^2)_{n}^3}
\bigg\{\sum_{i=1}^{n}\frac{q^{2j}}{[2j]^2}-\sum_{j=1}^{n+1}\frac{q^{2j-1}}{[2j-1]^2}\bigg\}.
\end{align*}
Letting $n\to \infty$ in this identity, we find \eqref{eq:wei-h}.
\end{proof}

Finally, we begin to prove Theorem \ref{thm-f}.

\begin{proof}[{\bf{Proof of Theorem \ref{thm-f}}}]
Setting $(a,b,c,d,e)=(q^{1/2},q^{1/2},x,q/x,q^{-n})$ in
\eqref{basic-a} and calculating the series on the right-hand side by
the $q$-analogue of Dougall's $_5F_4$ summation formula  (cf.
\cite[Appendix II. 21]{Gasper}):
\begin{align*}
{_6}\phi_{5}\left[\begin{array}{c}
a,qa^{\frac{1}{2}},-qa^{\frac{1}{2}},b,c,q^{-n}\\[3pt]
a^{\frac{1}{2}},-a^{\frac{1}{2}},aq/b,aq/c,aq^{n+1}
\end{array};\, q, \frac{aq^{n+1}}{bc} \right]
=\frac{(aq,aq/bc;q)_{n}}{(aq/b,aq/c;q)_{n}},
\end{align*}
 there is
\begin{align*}
&\sum_{k=0}^{n}\frac{(q^{\frac{1}{2}};q)_{k}^2(q^{1+n},q^{-n};q)_{k}}{(q,q^{\frac{1}{2}+n},q^{\frac{3}{2}+n};q)_{2k}}
\frac{(x;q)_{k}^2(q/x;q)_{k}^2(xq^{\frac{1}{2}+n},q^{\frac{3}{2}+n}/x;q)_{k}}{(xq^{\frac{1}{2}},q^{\frac{3}{2}}/x;q)_{2k}}
\\[3pt]
&\quad\times(-1)^kq^{\frac{5k^2+k}{2}+kn}U_k(x,n;q)=\frac{(q^{\frac{1}{2}},q^{\frac{3}{2}};q)_{n}}{(xq^{\frac{1}{2}},q^{\frac{3}{2}}/x;q)_{n}},
\end{align*}
where
\begin{align*}
&U_k(x,n;q)=\frac{1-q^{\frac{1}{2}+3k}}{1-q^{\frac{1}{2}}}+\frac{q^{\frac{1}{2}+2k+n}(1-q^{k-n})(1-xq^{k})(1-q^{1+k}/x)(1-xq^{\frac{3}{2}+3k+n})}
{(1-q^{\frac{1}{2}})(1+q^{\frac{1}{2}+k})(1-q^{\frac{1}{2}+2k+n})(1-q^{\frac{3}{2}+2k+n})(1-xq^{\frac{1}{2}+2k})}
\end{align*}
\begin{align*}
&\:+\frac{q^{2+4k+n}(1-q^{1+k+n})(1-q^{k-n})(1-xq^{k})^2(1-q^{1+k}/x)(1-xq^{\frac{1}{2}+k+n})}
{x(1-q^{\frac{1}{2}})(1+q^{\frac{1}{2}+k})(1-q^{\frac{1}{2}+2k+n})(1-q^{\frac{3}{2}+2k+n})^2(1-xq^{\frac{1}{2}+2k})(1-q^{\frac{3}{2}+2k}/x)}.
\end{align*}
Via the operator $\mathcal{D}_{x}$ and the last equation, it is
clear that
\begin{align}
&\sum_{k=0}^{n}\frac{(q^{\frac{1}{2}};q)_{k}^2(q^{1+n},q^{-n};q)_{k}}{(q,q^{\frac{1}{2}+n},q^{\frac{3}{2}+n};q)_{2k}}
\frac{(x;q)_{k}^2(q/x;q)_{k}^2(xq^{\frac{1}{2}+n},q^{\frac{3}{2}+n}/x;q)_{k}}{(xq^{\frac{1}{2}},q^{\frac{3}{2}}/x;q)_{2k}}
\notag\\[2pt]
&\quad\times(-1)^kq^{\frac{5k^2+k}{2}+kn}U_k(x,n;q)V_k(x,n;q)
\notag\\[2pt]
&+\sum_{k=0}^{n}\frac{(q^{\frac{1}{2}};q)_{k}^2(q^{1+n},q^{-n};q)_{k}}{(q,q^{\frac{1}{2}+n},q^{\frac{3}{2}+n};q)_{2k}}
\frac{(x;q)_{k}^2(q/x;q)_{k}^2(xq^{\frac{1}{2}+n},q^{\frac{3}{2}+n}/x;q)_{k}}{(xq^{\frac{1}{2}},q^{\frac{3}{2}}/x;q)_{2k}}
\notag\\[2pt]
&\qquad\times(-1)^kq^{\frac{5k^2+k}{2}+kn}\mathcal{D}_{x}U_k(x,n;q)
\notag\\[2pt]
&=\frac{(q^{\frac{1}{2}},q^{\frac{3}{2}};q)_{n}}{(xq^{\frac{1}{2}},q^{\frac{3}{2}}/x;q)_{n}}
\bigg\{\sum_{j=1}^n\frac{q^{j-\frac{1}{2}}}{1-xq^{j-\frac{1}{2}}}-\sum_{j=1}^n\frac{q^{j+\frac{1}{2}}/x^2}{1-q^{j+\frac{1}{2}}/x}\bigg\},
\label{eq:wei-aaa}
\end{align}
where
\begin{align*}
V_k(x,n;q)&=2\sum_{i=1}^k\frac{q^{i}/x^2}{1-q^{i}/x}-2\sum_{i=1}^k\frac{q^{i-1}}{1-xq^{i-1}}+\sum_{i=1}^{2k}\frac{q^{i-\frac{1}{2}}}{1-xq^{i-\frac{1}{2}}}
\\[2pt]
&\quad-\sum_{i=1}^{2k}\frac{q^{i+\frac{1}{2}}/x^2}{1-q^{i+\frac{1}{2}}/x}+\sum_{i=1}^{k}\frac{q^{i+\frac{1}{2}+n}/x^2}{1-q^{i+\frac{1}{2}+n}/x}
-\sum_{i=1}^{k}\frac{q^{i-\frac{1}{2}+n}}{1-xq^{i-\frac{1}{2}+n}}.
\end{align*}
Dividing both sides of \eqref{eq:wei-aaa} by $1-q/x^2$ and then
letting $(x,q,n)\to(q,q^2,\infty)$, there holds
\begin{align}
&\sum_{k=0}^{\infty}q^{6k^2}\frac{(q;q^2)_k^6}{(q^2;q^2)_{2k}^3}\bigg\{\sum_{i=1}^{2k}\frac{q^{2i}}{[2i]^2}
-2\sum_{i=1}^{k}\frac{q^{2i-1}}{[2i-1]^2}\bigg\}
\notag\\[2mm]
&\quad\times\frac{1+2q^{1+2k}-q^{1+6k}(2+2q^2+q^{1+2k}+q^{3+2k}-3q^{3+6k})}{(1-q)(1-q^{2+4k})(1+q^{1+2k})^2}
\notag\\[2mm]
&\:+\sum_{k=0}^{\infty}q^{6k^+8k+2}\frac{(q;q^2)_k^6}{(q^2;q^2)_{2k}^3}
\frac{(1-q)(1+q^{1+2k}+q^{2+4k})^2}{(1-q^{2+4k})(1+q^{1+2k})^4}
\notag\\[2mm]
&\:=\frac{(q,q^{3};q^2)_{\infty}}{(q^{2};q^{2})_{\infty}^2}\sum_{j=1}^{\infty}\frac{q^{2j}}{[2j]^2}.
 \label{eq:wei-bbb}
\end{align}

The $(a,b,c)=(x,q/x,q)$ case of \eqref{basic-d} can be stated as
\begin{align}
&\sum_{k=0}^{\infty}\frac{(x,q/x;q)_{k}^3}{(q^2;q)_{2k}^3}
 q^{3k^2}\bigg\{(1-q^{1+2k})^3(1-xq^{3k})-\frac{q^{1+3k}}{x}(1-xq^k)^3(1-q^{2+3k}/x)\bigg\}
 \notag\\[2pt]
&\:=(1-q)\frac{(x,q/x;q)_{\infty}}{(q^2;q)_{\infty}^2}.
\label{series-b}
\end{align}
 When $0<x<1$, it is obvious that the series on the
left-hand side of \eqref{series-b} is uniformly convergent. Through
the operator $\mathcal{D}_{x}$ and \eqref{series-b}, it is not
difficult to see that
\begin{align}
&3\sum_{k=0}^{\infty}\frac{(x,q/x;q)_{k}^3}{(q^2;q)_{2k}^3}
 q^{3k^2}\bigg\{(1-q^{1+2k})^3(1-xq^{3k})-\frac{q^{1+3k}}{x}(1-xq^k)^3(1-q^{2+3k}/x)\bigg\}
 \notag\\[2pt]
&\quad\times\bigg\{\sum_{i=1}^k\frac{q^{i}/x^2}{1-q^{i}/x}-\sum_{i=1}^k\frac{q^{i-1}}{1-xq^{i-1}}\bigg\}
\notag\\[2pt]
&\:+\sum_{k=0}^{\infty}\frac{(x,q/x;q)_{k}^3}{(q^2;q)_{2k}^3}
q^{3k^2+3k}\frac{(q-x^2)(x+3xq^{2+4k}-2x^2q^{1+3k}-2q^{2+3k})}{x^3}
 \notag\\[2pt]
&\:=(1-q)\frac{(x,q/x;q)_{\infty}}{(q^2;q)_{\infty}^2}
\bigg\{\sum_{j=1}^{\infty}\frac{q^{j}/x^2}{1-q^{j}/x}-\sum_{j=1}^{\infty}\frac{q^{j-1}}{1-xq^{j-1}}\bigg\}.
 \label{eq:wei-ccc}
\end{align}
Dividing both sides of \eqref{eq:wei-ccc} by $1-q/x^2$ and then
letting $(x,q)\to(q,q^2)$, we can verify that
\begin{align}
&3\sum_{k=0}^{\infty}q^{6k^2}\frac{(q;q^2)_{k}^6}{(q^4;q^2)_{2k}^3}
 \frac{(1-q^{2+4k})^3(1-q^{1+6k})-q^{1+6k}(1-q^{1+2k})^3(1-q^{3+6k})}{(1-q)^4}
 \notag\\[2pt]
&\quad\times\sum_{i=1}^{k}\frac{q^{2i-1}}{[2i-1]^2}
+\sum_{k=0}^{\infty}q^{6k^2+6k}\frac{(q;q^2)_{k}^6}{(q^4;q^2)_{2k}^3}
\frac{q-4q^{4+6k}+3q^{5+8k}}{(1-q)^2}
 \notag\\[2pt]
&\:\:=(1+q)^2\frac{(q^3;q^2)_{\infty}^2}{(q^2,q^4;q^2)_{\infty}}
\sum_{j=1}^{\infty}\frac{q^{2j-1}}{[2j-1]^2}.
 \label{eq:wei-ddd}
\end{align}
By means of the linear combination of \eqref{eq:wei-bbb} and
\eqref{eq:wei-ddd} multiplied, respectively, by $(-64)$ and $3$, we
are led to \eqref{eq:wei-i}.
\end{proof}

%%%%%%%%%%%%%%%%%%%%%%%%%%%%%%%%%%%%%%%%%%%%%%%%%%%%%%%%%%%%%%%%%%%%%%%%%%%%%%%%%%%%%%%%%%%%%%%%%%%%%%%%%%%%%%%%%%%%%%%%%%%%%%%%%%%%%%%%%%%%%%%%%%%%%%%%%%%%%%%%%%%%%%%%%%%%%%%%%%%%%%%%%%%%%%%%%%%%%%%%%%

\end{document}